\documentclass[12pt,a4paper,twoside,final,notitlepage, leqno]{article}
\usepackage[english]{babel}
\usepackage[T1]{fontenc}  
\usepackage{epsfig, graphicx, subfigure, amssymb}
\usepackage{amsmath,amsthm,epsfig,amsfonts,bbm}  

\newcommand{\wt}{\widetilde}
\newcommand{\ds}{\displaystyle}

\newcommand{\pequationdeb}{$$ \left\{ \begin{minipage}[c]{130mm}}
\newcommand{\pequationfin}{\end{minipage}
                           \right. $$}

\newcommand{\beq}     {\begin{equation}}
\newcommand{\enq}     {\end{equation}}
\newcommand{\be}    {\begin{enumerate}}
\newcommand{\ee}    {\end{enumerate}}

\newcommand{\Bb}

\def\resume{\if@twocolumn
\section*{R\'esum\'e}
\else \small
\quotation{\bf \it R\'esum\'e \rule[1mm]{1.5mm}{0.2mm}\vspace{0pt}}
\fi}
\def\endresume{\if@twocolumn\else\endquotation\fi}
\def\abstract{\if@twocolumn
\noindent\section*{{\bf Abstract}}
\else \small
\quotation{\noindent \bf {Abstract.} \rule[1mm]{1.5mm}{0.2mm}\vspace{0pt}}
\fi}
\def\endabstract{\if@twocolumn\else\endquotation\fi}

\def\section*#1{}

\usepackage{fancyhdr}
\renewcommand{\headrulewidth}{0pt}
\begin{document} 

\fancypagestyle{plain}{ \fancyfoot{} \renewcommand{\footrulewidth}{0pt}}
\fancypagestyle{plain}{ \fancyhead{} \renewcommand{\headrulewidth}{0pt}}
\bibliographystyle{alpha}

\title{\bf \LARGE  On a superconvergent  lattice Boltzmann \\ boundary scheme   }

\author { {\large   Fran\c{c}ois Dubois  $^{a b}$,  Pierre Lallemand $^c$ and 
 Mahdi Tekitek  $^{a}$$^d$} \\~\\
{\it \small  $^a$   Numerical Analysis and Partial Differential Equations }   \\
 {\it  \small  Department of Mathematics, Paris Sud  University,  Orsay, France. } \\
{\it  \small $^b$ Conservatoire National des Arts et M\'etiers, LMSSC, Paris, France.} \\
{\it  \small $^c$  Beijing Computational Science Research Center, Beijing, China.}  \\
{\it  \small  $^d$ Present address: Departement of Mathematics, F.S.T., 
University El Manar, Tunis, Tunisia.}  \\
 {\rm  \small  francois.dubois@math.u-psud.fr,   pierre.lal@free.fr,
 mohamed-mahdi.tekitek@math.u-psud.fr} }  

\date { { \large   04 June 2009 } 
 \footnote {\rm  \small $\,\,$ Contribution presented at the fifth  ICMMES Conference 
(Amsterdam, 16-20 June 2008) and published in   {\it Computers and Mathematics with 
Applications},  volume 59, number 7,   pages  2141-2149, 
   april~2010,  doi:10.1016/j.camwa.2009.08.055. Edition  04 may 2014.  }} 

\maketitle

\begin{abstract}
\noindent In a seminal paper \cite{bibl:GA94} Ginzburg and Adler
 analyzed the bounce-back boundary conditions for the lattice Boltzmann
scheme   and showed that it could
 be made exact to second order 
  for  the Poiseuille flow  if some expressions depending upon the parameters of
 the method were satisfied,  thus defining so-called ``magic
 parameters''.   Using the Taylor expansion   method that one
 of us developed, we analyze a series of simple situations (1D and 2D) for
 diffusion and for linear fluid problems  using  bounce-back and 
``anti bounce-back'' numerical boundary conditions. 
The result is that ``magic parameters''  depend upon the
 detailed choice of the moments and of their equilibrium values. They may also
 depend upon the way the flow is driven.\\
{\bf Keywords:} 
Lattice Boltzmann scheme, boundary conditions, Taylor expansion method.   \\ 
   {\bf AMS classification}: 65-05, 65Q99, 82C20.   
\end{abstract}  

\newpage 
\fancyfoot[C]{\oldstylenums{\thepage}}
 
\bigskip  \noindent {\bf \large 1 $\,\,$ Introduction }  

\noindent  
The theoretical analysis of the lattice Boltzmann scheme \cite{bibl:FHH87, bibl:HJ89, bibl:HSB89,
bibl:qian, bibl:dh92, bibl:KGSB98,  bibl:La00}
is an active subject of research.
Recall that the   method  was  first analyzed by d'Humi\`eres
\cite{bibl:dh92}  with a Chapman Enskog expansion coming from statistical physics;  
we also refer to Asinari and Ohwada~\cite{bibl:AO09} for a method of analysis 
 based on the Grad moment system.
A fruitful idea followed by 
Junk {\it et al}~\cite{bibl:JKL05} and~\cite{bibl:Du07,bibl:Du08}
is to use the so-called equivalent equation method derived independently by 
Lerat-Peyret~\cite{bibl:LP74} and Warming and Hyett~\cite{bibl:WH74}
(see also \cite{bibl:VR99}). An infinitesimal  parameter is introduced and 
the finite differences operators are expanded 
into a family of equivalent partial differential equations. 
The main goal of this study is to use the Taylor expansion 
method~\cite{bibl:Du07,bibl:Du08} in order to increase the accuracy
of   boundary conditions for simple problems with analytical
solutions. 
We first consider a one-dimensional (1D) diffusion problem
and study the influence of   the definition of the moments and of
their equilibrium value.
We then consider a two-dimensional (2D) Poiseuille flow using 
several ways to enforce a pressure gradient. 

We consider regular lattices parametrized by a space step 
 $\Delta{x}$. We introduce a time step  $\Delta{t}$ and adopt the 
 ``acoustic'' scaling: the ratio 
 $\lambda \equiv \frac{\Delta{x}}{\Delta{t}}$ 
is a  {\bf fixed} reference velocity for each study. 
As a consequence, the parameters  $\Delta{x}$ and  $\Delta{t}$ are equivalent
 infinitesimals. 
Note that as this work is devoted to boundaries, we shall use a particular way to
test the accuracy of a numerical scheme as will be discussed later.

\bigskip  \noindent {\bf \large 2 $\,\,$ Diffusion problem in one space dimension} 

\noindent   
We consider the classical Lattice Boltzmann model D$1$Q$3$ with three 
discrete velocities and  one conservation law to model diffusion problems.  
%
We choose the velocities $v_{i}$ $(0 \leq i \leq 2)$  such that 
$\, v_0 = 0, \, v_1 = \lambda , \, v_2=- \lambda$. 
At each mesh point, there are three functions $\{f_j\}$
that can be interpreted as populations of fictitious particles.
These populations evolve according to the lattice Boltzmann scheme which 
we write as in \cite{bibl:Du07}:
\begin{equation}
f_{j}(x,\, t+\Delta{t})=f^{*}_{j}(x - v_{j} \Delta{t}, \, t) 
\,, \qquad 0\leq j \leq 2,  \label{dubois}
\end{equation}
where the superscript $*$ denotes post-collision quantities and $x$ a vertex of the
lattice. Therefore
during each  time increment $\Delta{t}$ there are two fundamental steps: advection
and collision.
 The {\bf advection}  step describes the motion of a particle which has undergone
collision at node $x- v_{j} \Delta{t}$ 
and goes to the $j$th  neighbouring node.
Following d'Humi\`eres~\cite{bibl:dh92}, the  {\bf collision} step is defined in
the space of moments. For D$1$Q$3$   three  moments $\{ m_{\ell} \}$ 
 are obtained by a linear transformation of vectors $f_{j}$:
\begin{equation}
m_{0} \,=\,  f_{0}+f_{1}+f_{2} \equiv  \rho \,\, ( \textrm {density})    , \quad  
m_{1}   \,=\, \lambda(f_{1}-f_{2}) , \quad  
m_{2}   \,=\,  \frac{\lambda^{2}}{2}(f_{1}+f_{2}) .   
\label{moments-d1q3}
\end{equation}
In consequence,
we introduce a matrix of moments $M$ to represent moments like (\ref{moments-d1q3});
it  takes the form
\begin{equation}
{{M}}=\left ( \begin{array}{ccc} 
1  &\quad 1 &\quad 1 \\  
0  &\lambda & -\lambda\\
0  &\frac{\lambda^2}{2} & \frac{\lambda^2}{2}
\end{array}  \right) 
\label{matrixM}
\end{equation}
and the relations  (\ref{moments-d1q3}) can be simply written as $ \, m = M \, f$. 
To simulate diffusion problems, we conserve only the density moment $\rho$ in
the collision step and obtain one macroscopic scalar equation. The other quantities
(non-conserved moments) are assumed to relax
towards equilibrium values ($m_{1}^{eq}$, $m_{2}^{eq}$)
following:
\begin{equation}
m_{\ell}^{*}=(1-s_{\ell})\,  m_{\ell} + s_{\ell} \, m_{\ell}^{eq},  \quad 1\leq \ell
 \leq 2,
\label{veq1}
\end{equation}
where $s_{\ell}$ ($ 0 <  s_{\ell} < 2 $, for $\ell =1,\, 2$) are relaxation rates,  not
necessarily equal to a single value as in the BGK case~\cite{bibl:qian}.
The equilibrium values $m_{\ell}^{eq}$ of the non conserved moments in 
equation $(\ref{veq1})$ determine the macroscopic behavior of the
scheme. 
Indeed with the
following choice of equilibrium values (neglecting non-linear contributions):
\begin{equation}
m_{1}^{eq} \,=\, 0 \,, \quad  m_{2}^{eq} = \zeta \frac{\lambda^{2}}{2} \rho
\label{meqd1q3}
\end{equation}
and using the Taylor expansion method  we find (see {\it e.g.} \cite{bibl:Du09}) that 
 the  equivalent partial differential equation of the numerical scheme 
up to order three in $\Delta x$  is a  diffusion equation: 
\begin{equation}
\ds \frac{\partial \rho}{\partial t} - 
\kappa \, \frac{\partial^{2} \rho}{\partial x^{2}} \,= \,  
\ds {\rm{O}}(\sigma_1 \, \Delta  x^{3}) \,.
\label{acou2}
\end{equation} 
%
The value of the diffusivity $\kappa$ is  given according to 
\begin{equation}
\kappa \,= \, \Delta t \, \lambda^{2} \, \sigma_{1} \, \zeta \,   
\label{diffusivite}
 \end{equation}
where $  \sigma_{\ell} \equiv  \frac{1}{s_{\ell}}-\frac{1}{2} \,, $ 
$\,  \ell = 1,\, 2.$ 

\smallskip 
Remark that  the thermal diffusivity  $\kappa$ is imposed by the Physics.
Moreover the scale velocity $ \lambda $ is fixed and the coefficient $ \zeta $
is also imposed. When we refine the mesh, 
the coefficient $  \sigma_{1} $
must be chosen in order to enforce  relation (\ref{diffusivite}).
In other terms, the product $ \,  \sigma_{1} \,   \Delta t \, $ 
must be maintained   constant. Then the right hand side of 
relation (\ref{acou2}) exhibits a {\bf second order
truncation error}  of the lattice Boltzmann scheme 
for a {\bf given}  thermal diffusivity   $\kappa$. 
Associated with  stability properties (see Junk and Yong  \cite {jy09b}),  
 convergence properties of lattice Boltzmann scheme can be established, 
as in  \cite{JY09a}.

\bigskip  \noindent {\bf \large 3 $\,\,$ Localization of a one-dimensional  boundary }

\noindent    
Let us introduce  a constant $c$ and  consider the following one-dimensional
Poisson  problem:
%
\begin{equation}
-{\bf{\rm{K}}}  \frac{{\rm d}^{2} \rho}{{\rm d}x^{2}}  \,=\, c  
\quad {\mbox{in}} \,  ]0,1[,  \qquad   \rho(0) \,=\,   \rho(1) \,=\, 0 \, .  
\label{pb1} 
\end{equation}   
%
\begin{figure}[!h]      
\centerline {\includegraphics [ width=.35 \textwidth ,  angle=0]{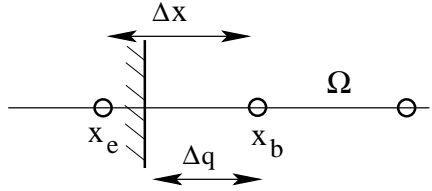}}
\caption{A boundary surface cutting the link between  a fluid node $x_{b}$ 
 and a fictitious outside  node $x_{e}\equiv x_{b}-\Delta x$. }
\label{diri}
\end{figure}
%
We take an ``anti bounce-back'' numerical boundary  condition at $x=0$:
\begin{equation}
f_{1}(x_{b},t+\Delta t)=-f_{2}(x_{e},t+\Delta t)=-f_{2}^{*}(x_{b},t) \,, 
\label{antibb-1D} 
\end{equation}
with $ \, x_{b} \,$ the fluid node and  $ \, x_{e} \,$ the external node as 
presented in   Figure $\ref{diri}$, and a similar condition for $x=1$. 
A  uniform body source ($\delta p$) is added to the Boltzmann scheme to model the right
hand side $c$ of equation (\ref{pb1}). So we can write the lattice Boltzmann scheme as follows: 
({\it i})~$ \, m = M \, f ,\,\,$ 
({\it ii})~$ \, \wt{m}_{0}=m_{0}+\frac{1}{2} \, \delta p ,\,\,$ 
({\it iii})~evaluate the other moments, $\,\,$ 
({\it iv})~relaxation (\ref{veq1}) of the other moments,  $\,\,$ 
({\it v})~$ \, \wt{m}_{0}=m_{0}+\frac{1}{2} \, \delta p , \,\,$  
({\it vi})~$ \, f = M^{-1} \, m , \,\,$  
({\it vii})~advection step  (\ref{dubois})  and boundary conditions.   
%
%
%
The exact solution of  problem $(\ref{pb1})$ is
elementary:   $\, u(x)= \frac{c \, x\, (1-x)}{2 \, {\rm K}} $.
We analyze the behavior of the discrete model to show whether it can
be tuned so that the location of the ``numerical boundary'' can be
fixed at mid-point as expected from ``anti bounce-back''. Thus we shall use as
criterion for accuracy the difference between the imposed boundary and the
``numerically determined'' boundary.

\smallskip 
From a theoretical point of view, we suppose that the discrete fields $ \, f_j(x,\, t) \,$
vary slowly in space and time in order to be able to use Taylor expansions.
We analyse the lattice Boltzmann scheme in terms of equivalent partial
differential  equations and formal developments.  
It is well known (see {\it e.g.}  Griffiths  and Sanz-Serna  \cite{GSS86} 
 or  Chang   \cite{bibl:Ch90})   that this method of analysis
 fails {\it a priori} to predict  
boundary effects properly if this hypothesis is not satisfied. 
We keep in mind this restriction in our numerical experiments. 
Nevertheless,  this elementary  tool can produce  nontrivial results, 
as we  will see hereafter. 

\smallskip 
We say in the following that  a boundary scheme 
(such as   (\ref{antibb-1D})  to fix the ideas) {\bf is of order $p$ at location $\Delta q$}
  relative to homogeneous  Dirichlet boundary condition 
(that are present  in (\ref{pb1})) if the
  numerical boundary condition implies 
\begin{equation}
\rho( x_b - \Delta q) \,=\, {\rm O}( \Delta x^p) 
\label{def-ordre-cl}  
\end{equation} 
for the continuous conserved field issued from the particle field $ \, f_j(x,\,t)\,$
according  to (\ref{moments-d1q3}). We have the following result: 

\smallskip \noindent  {\bf Proposition 1. \quad   Superconvergent   relation between parameters }

\noindent
For the D$1$Q$3$ lattice Boltzmann scheme  
 (\ref{dubois})    (\ref{moments-d1q3})   (\ref{veq1})     (\ref{meqd1q3}),  
the ``anti bounce-back'' numerical boundary  condition  (\ref{antibb-1D}) 
is of order $3$ at location $\Delta q = \frac{\Delta x}{2}\, $ 
relative to the homogeneous Dirichlet boundary condition of problem   (\ref{pb1}) if 
and only if   the following condition 
\begin{equation} 
\sigma_{1} \, \sigma_{2} \,=\, \frac{1}{8} 
\label{magic-irina}  
\end{equation} 
%
is satisfied.   

\smallskip \noindent  
Relation (\ref{magic-irina}) defines superconvergent  parameters $\sigma_1$ and  $\sigma_2$. 
Recall that they have been called ``magic'' by I.~Ginzbourg and  P.M.~Alder~\cite{bibl:GA94}.
%

\smallskip \noindent   {\bf Proof of Proposition 1.}

\noindent  
 We have introduced in  
\cite{bibl:Du07, bibl:Du08}   the ``tensor of momentum velocities'' 
and the so-called ``defects of conservation'' which are defined respectively by
\begin{equation*}  
\Lambda_{k\alpha}^{\ell} \,\equiv \, \sum_j M_{k j}  
\, M_{\alpha j} \,  (M^{-1})_{j \ell} \, 
\end{equation*} 
%
%
\begin{equation}
\theta_{k}   \, \equiv \, \partial_{t} m_{k}^{eq}+ \Lambda_{k\alpha}^{\ell} \, 
 \partial_{\alpha} m_{\ell}^{eq}, \qquad k \geq 1.  \label{defect}
\end{equation}
For the  D$1$Q$3$ lattice Boltzmann scheme applied to diffusion problem like 
(\ref{pb1}),
\begin{equation}
\theta_{1}  \, = \,  \zeta \lambda^{2} \frac{\partial \rho}{\partial {x}} \,, \qquad 
\theta_{2}    \, = \,  \zeta \, \frac{\lambda^{2}}{2} \frac{\partial \rho}{\partial {t}} \, .
\label{t1}  
\end{equation} 
Then we obtain the following development of non-equilibrium moments at third
order  (as described in \cite{bibl:Du09}):
\begin{equation}
m_{k}^{*}=m_{k}^{eq}+\Delta t \, 
 \big( \frac{1}{2}-\sigma_{k}  \big) \, \left[\theta_{k}-\Delta t
  (\sigma_{k}\partial_{t} \theta_{k}+\sigma_{\ell}
  \Lambda_{k\alpha}^{\ell}\partial_{\alpha}\theta_{\ell} )\right]   +{\rm{O}}(\Delta x^{3}),
 \quad k \geq 1.
\label{order2}
\end{equation}
Thus for $k=1$:
\begin{equation*}
m_{1}^{*} = \Delta t \, 
 \Big(\frac{1}{2}-\sigma_{1} \Big) 
\left[\theta_{1}-\Delta t( \sigma_{1}\partial_{t}\theta_{1}+\sigma_{1}
 \Lambda_{11}^{1}\partial_{x}\theta_{1}+\sigma_{2}
 \Lambda_{11}^{2}\partial_{x}\theta_{2})\right]+{\rm{O}}(\Delta x^{3}) \, .  
\end{equation*} 
Using $(\ref{t1})$ and  $\partial_{t} \theta_{1}={\rm{O}}(\Delta t^{2})$, 
$\Lambda_{11}^{1}=0$
the above equation becomes:
\begin{equation}
m_{1}^{*}=\Delta t \lambda^{2} \left(\frac{1}{2}-\sigma_{1}\right) \zeta
 \frac{\partial {\rho}}{\partial x}+{\rm{O}} (\Delta x^{3}). \label{qq}
\end{equation}
For $k=2$, we use expression $(\ref{meqd1q3})$ of $m_{2}^{eq}$,
together with $\theta_{2}=0$ and $\Lambda_{21}^{1}=\frac{\lambda^{2}}{2}$
to obtain from  equation ($\ref{order2}$):
\begin{equation}
m_{2}^{*}=\lambda^{2} \, \frac{\zeta}{2} \,  \rho \,- \, \Delta t^{2} \, \lambda^{4}
\,  \frac{\zeta}{2} \, \sigma_{1} \,   \big( \frac{1}{2}-\sigma_{2}  \big) 
 \frac{\partial^2 \rho}{\partial x^{2}}+{\rm{O}} (\Delta x^{3}) \label{ee}
\end{equation}
Using the inverse moment matrix $M^{-1}$, the post-collision $f$ are given by:
\begin{eqnarray}
f_{1}^{*}=\frac{1}{2\lambda^{2}}\left[2 m_{2}^{*}+\lambda m_{1}^{*} \right], \quad
f_{2}^{*}=\frac{1}{2\lambda^{2}}\left[2 m_{2}^{*}-\lambda m_{1}^{*} \right].
\end{eqnarray}
At the boundary, due to (\ref{antibb-1D}) and  (\ref{dubois}),  we consider the following quantity:
\begin{equation} 
f_{1}^{*}(x_{e})+f_{2}^{*}(x_{b})=\frac{1}{2\lambda^{2}} 
\left[2(m_{2}^{*}(x_{e})+m_{2}^{*}(x_{b}))+\lambda(m_{1}^{*}(x_{e})-m_{1}^{*}(x_{b}))\right].
\label{diff}
\end{equation}
Using relations $(\ref{ee})$ and $(\ref{qq})$ we obtain respectively:
\begin{equation} \label{in1}  \left\{ \begin{array}{c} \displaystyle  
m_{2}^{*}(x_{e})+m_{2}^{*}(x_{b})   
\,=\,   \lambda^{2} \, \frac{\zeta}{2} \, \left[ \rho(x_{e})+\rho(x_{b}) \right] 
 \qquad  \qquad \qquad    \qquad    \\   [1mm] 
 \displaystyle   \qquad \qquad  
 - \Delta t^{2} \, \lambda^{4} \, \frac{\zeta}{2}\sigma_{1} \,   
\big(\frac{1}{2} \, - \, \sigma_{2}\big) \, \Big[ \frac{\partial^{2}\rho}{\partial
 x^{2}}(x_{e}) \, + \frac{\partial^{2}\rho}{\partial x^{2}}(x_{b}) \Big] 
+{\rm{O}}  (\Delta x^{3})  
\end{array}  \right.  \end{equation}  
\begin{equation} 
m_{1}^{*}(x_{e})-m_{1}^{*}(x_{b}) \,=\,     \Delta t \, \lambda^{2} \, \zeta \, 
 (\frac{1}{2}-\sigma_{1})   \, \left[ \frac{\partial \rho}{\partial
 x}(x_{e})-\frac{\partial \rho}{\partial x} (x_{b}) \right] 
+{\rm{O}} (\Delta  x^{3}) \label{in2}.
\end{equation} 
With the help of classical Taylor expansion we have, with the notation 
$ \, x_i \equiv {1\over2} (x_b+x_e) $:
\begin{eqnarray}
\rho(x_{e}) + \rho(x_{b}) &=& 2\rho(x_{i})+\frac{\Delta
 x^{2}}{4}\, \frac{\partial^{2} \rho}{\partial x^{2}}(x_{i})+{\rm{O}}
 (\Delta x^{3}),\label{in4}\\
\frac{\partial \rho}{\partial  x}(x_{e}) -  \frac{\partial \rho}{\partial x} (x_{b})  
&=& -\Delta x \, 
 \frac{\partial^{2} \rho}{\partial x^{2}}(x_{i})+{\rm{O}} (\Delta x^{3}),\label{in3}\\
 \frac{\partial^{2}\rho}{\partial
  x^{2}}(x_{e}) + \frac{\partial^{2}\rho}{\partial
  x^{2}}(x_{b})  &=&  2 \, \frac{\partial^{2}\rho}{\partial x^{2}}(x_{i})+{\rm{O}}
  (\Delta x^{2}).\label{in5} 
\end{eqnarray}
Considering equation $(\ref{diff})$, together with ($\ref{in1}$),  ($\ref{in2}$) 
and taking into account relations 
($\ref{in4}$),  ($\ref{in3}$) and ($\ref{in5}$) we obtain:
\begin{equation}
f_{1}^{*}(x_{e})+f_{2}^{*}(x_{b}) \,=\, \zeta \,\rho(x_{i}) +   \zeta \, \Delta x^{2}  
  \left( \sigma_{1}\sigma_{2}-\frac{1}{8} \right)  \, 
\frac{\partial^{2}\rho}{\partial x^{2}}(x_{i})+{\rm{O}}
 (\Delta x^{3}).
\label{dvpt-d1q3}  
\end{equation}
Due to the simple fundamental expression (\ref{dubois})  of a lattice Boltzmann
scheme,   the left hand side of 
(\ref{dvpt-d1q3})  is identically  null 
when the numerical boundary condition (\ref{antibb-1D})  occurs. 
 Due to the relation  
$x_i =  {1\over2} (x_b+x_e) $, the condition 
(\ref{def-ordre-cl}) is satisfied with   
$ \, \Delta q =\frac{\Delta x}{2}$ and $p=3$ if and only if    
$ \, \sigma_{1}\sigma_{2}=\frac{1}{8} .  $  
 \hfill $\Box$ 

\bigskip   \noindent  $\bullet$ \quad 
Let us now consider the effect of using a {\bf different}  moment matrix for the D$1$Q$3$ case:
\begin{equation}
{{M}}=\left ( \begin{array}{ccc} 
 1  &\quad 1 &\quad 1 \\  
 0  &\lambda & -\lambda\\
 -2 \lambda^{2} &\lambda^{2} & \lambda^{2}
\end{array}  \right).
\label{matrixMM}
\end{equation}
obtained from (\ref{matrixM}) by a Gram-Schmidt orthogonalization algorithm as usual with
the lattice Boltzmann scheme (see {\it e.g.} \cite{bibl:La00}).  
The moments at equilibrium are now  given by $ \,  m_{1}^{eq} = 0 \, $ and 
$ \,  m_{2}^{eq} = \lambda^{2} \,  \wt{\zeta} \, \rho. $ 
We remark that the matrix of moments ($\ref{matrixMM}$) 
leads to an equivalent macroscopic conservation law of type 
(\ref{acou2}) with 
a diffusivity  $ \kappa $ which is now given by  
 $  \kappa  = \Delta t \, \lambda^{2} \, \sigma_{1} \, \frac{2+\wt{\zeta}}{3}$.

\smallskip  \noindent  {\bf Proposition 2. \quad   Third order at the boundary }

\noindent
The  D$1$Q$3$ lattice Boltzmann scheme  
 (\ref{dubois})    (\ref{matrixMM})   (\ref{veq1})     (\ref{meqd1q3}) 
associated to the ``anti bounce-back'' numerical boundary  condition  (\ref{antibb-1D}) 
is of order $3$  at location $\Delta q = \frac{\Delta x}{2}\, $ 
for  the homogeneous Dirichlet boundary condition of problem   (\ref{pb1}) if  
and only if    $ \, \sigma_{1}\sigma_{2}=\frac{3}{8}. \, $ 
%

\smallskip \noindent 
Remark that the superconvergent parameters satisfying the relation  (\ref{magic-irina}) 
 emerging from  Proposition 1  with the choice of transformation matrix  $M$  
  given by (\ref{matrixM})  
 are  {\bf different} from those obtained  in the case with  matrix 
(\ref{matrixMM}). 
%

\smallskip \newpage \noindent   {\bf Proof of Proposition 2.}

\noindent   
For this model we have the following Taylor development of non-conserved moments up to
 order $2$ on $\Delta x$:
\begin{eqnarray*}
m_{1}^{*}&=&\Delta t \, \lambda^{2} \, \big(\frac{1}{2}-\sigma_{1}\big) \, 
 \Big(\frac{2+\wt{\zeta}}{3}\Big) \, 
 \frac{\partial \rho}{\partial x} \,+\, {\rm{O}} (\Delta x^{3}),\\
m_{2}^{*}&=&\lambda^{2} \, \wt{\zeta} \, \rho - \Delta t^{2} \, \lambda^{4} \, 
 \sigma_{1}  \, \big(\frac{1}{2}-\sigma_{2}\big) \, 
 \Big(\frac{2+\wt{\zeta}}{3}\Big) 
\,  \frac{\partial^{2} \rho}{\partial x^{2}}  \,+\, {\rm{O}} (\Delta x^{3}).
\end{eqnarray*} 
With the help of the matrix moments ($\ref{matrixMM}$) we have:
\begin{equation*}
f_{1}^{*} = \frac{1}{3} \, \rho \,+\, 
\frac{1}{6\lambda^{2}} \, \big[m_{2}^{*}+3\lambda m_{1}^{*} \big], \qquad
f_{2}^{*} =   \frac{1}{3} \, \rho \, + \, \frac{1}{6\lambda^{2}} \, 
\big[m_{2}^{*}-3\lambda m_{1}^{*} \big].
\end{equation*}
As $f_{2}(x_{b})$ is internal to the domain we add to $\rho(x_{b})$ a 
 body source $\delta \rho  = -\frac{2+\wt{\zeta}}{3} \, 
 \frac{\partial^{2} \rho}{\partial x^{2}}$.  
Now by using the same method as in the proof of Proposition 1, we obtain:
\begin{equation*}
f_{1}^{*}(x_{e})+f_{2}^{*}(x_{b}) \,=\, \frac{2+\wt{\zeta}}{3} \rho(x_{i})
+\Delta x^{2} \, \big( 8 \sigma_{1} \sigma_{2} -3 \big) 
\, \frac{2+\wt{\zeta}}{72} \, \frac{\partial^{2}\rho}{\partial x^{2}}(x_{i})+{\rm{O}}
 (\Delta x^{3}).
\end{equation*}
The conclusion is a direct consequence of the above calculus. 
\hfill $\Box$ 

\begin{figure}[!h]      
\centerline {\includegraphics [ width=.50 \textwidth ,  angle=-90]{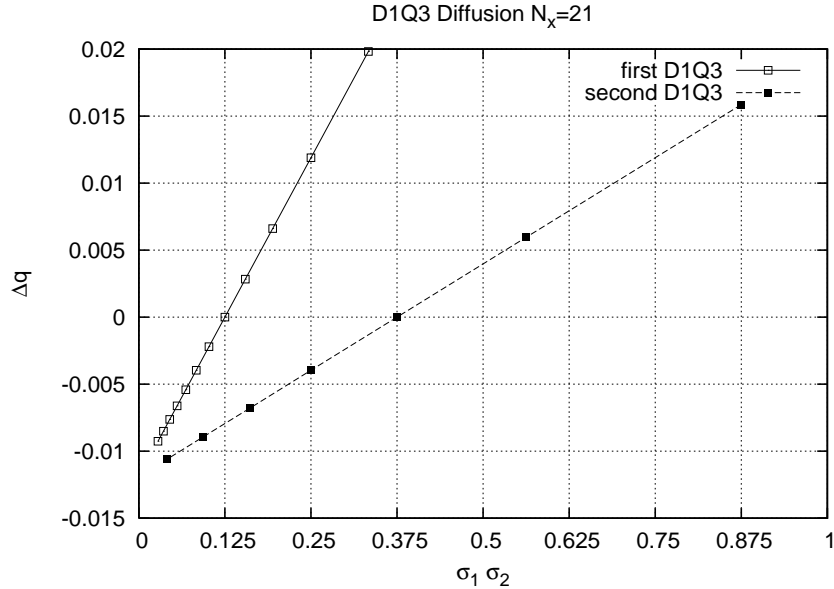}}
 \caption{ ``Experimental numerical  location'' of the solid wall  $\Delta q$  {\it versus} 
 $\sigma_{1}\sigma_{2}$. 
 The first D$1$Q$3$ model with $ \Box$ symbol: superconvergent parameters such that 
 $\sigma_{1}  \, \sigma_{2}=\frac{1}{8}$.  
The second D$1$Q$3$ model with  $\blacksquare$  symbol:  superconvergent parameters
satisfying    $\sigma_{1} \, \sigma_{2}=\frac{3}{8}$.}  
\label{compd1q3}
\end{figure}
%

\smallskip \noindent 
To illustrate the preceeding discussion, we perform a numerical simulation of the
two lattice Boltzmann  models and analyze (after a suitable number of iterations
to reach steady state) the ``Poiseuille'' parabolic profile. 
We measure the numerical  error in terms 
of a   precise location of the boundary for Dirichlet type boundary condition. 
We follow  a method proposed by   Ginzburg and  d'Humi\`eres~\cite{bibl:GD95}:
from the numerical discrete field $ \, u_{LB}(j \, \Delta x) \,$ we determine by least squares
a parabola that fit at best the data. Then  we calculate where this approximation of the 
 numerical  solution $u_{LB}$ is equal to zero. We interpret this location as 
 the ``experimental numerical  location''  of the solid wall. We find {\bf experimentally } 
  that the extrapolated  location of the Dirichlet boundary condition 
is located between $x_{b}$ and $x_{e}$ 
and this exact solid wall location  is parametrized under the form  
$x_{b}-\Delta q$, with  $0 \leq \Delta q \leq \Delta x$.
The results obtained for
several values of $\sigma_1$ and $\sigma_2$ are shown in Figure \ref{compd1q3}  to depend only upon 
the product $\sigma_1 \sigma_2$ and  go through 0 respectively for 
$\frac{1}{8}$ or  $\frac{3}{8}$, in {\bf  complete coherence}
 with the Taylor expansion method developed in Propositions 1 and 2. 
%

\bigskip  \noindent {\bf \large 4 $\,\,$ The two-dimensional Poiseuille flow}

\noindent   
We consider here the classical D$2$Q$9$ model (see {\it e.g.} 
\cite{bibl:La00}). We  study a Poiseuille flow (in linear regime),
first  with an imposed uniform body 
force and periodic boundary condition at the inlet and oulet of the
channel. Then we consider the same flow 
with an imposed difference of pressure between inlet and outlet.
%
The evolution of the lattice Boltzmann scheme is given by equation (\ref{dubois}). 
The corresponding moments have an explicit physical significance: 
$m_0\equiv \rho$ is the density, $m_{1}\equiv j_x$ and
$m_{2} \equiv j_{y}$ are $x$ and $y$ components of momentum, $m_3$ is the
energy, $m_4$ is related to square energy, $m_5, m_6$ are $x$ and $y$
components of heat flux and $m_7$, $m_8$ are diagonal stress and off-diagonal
stress.  A Gram-Schmidt orthogonalization method is also used and the matrix of moments is
exactly that used in \cite{bibl:Du07, bibl:La00}. 
The collision is described in the moments space as:
\begin{equation}
m_{\ell}^{*}=(1-s_{\ell}) \, m_{\ell} + s_{\ell} \, m_{\ell}^{eq},  \quad 3\leq \ell
 \leq 8,
\label{moments-hors-eq}
\end{equation}
where the equilibrium values $m_{k}^{eq}$ are given by:
\begin{equation}
m_{3}^{eq}=\alpha \, \rho \,,\,
 m_{4}^{eq}=\beta \, \rho  \,,\,
m_{5}^{eq}=-\frac{j_{x}}{\lambda} \,,\, 
 m_{6}^{eq}=-\frac{j_{y}}{\lambda} \,,\, 
 m_{7}^{eq}=0  \,,\,   m_{8}^{eq}=0 \, . 
\label{eq-d2q9}
\end{equation}
%
%
\begin{figure}[!h]      
\centerline { {\includegraphics [ width=.4  \textwidth  ]{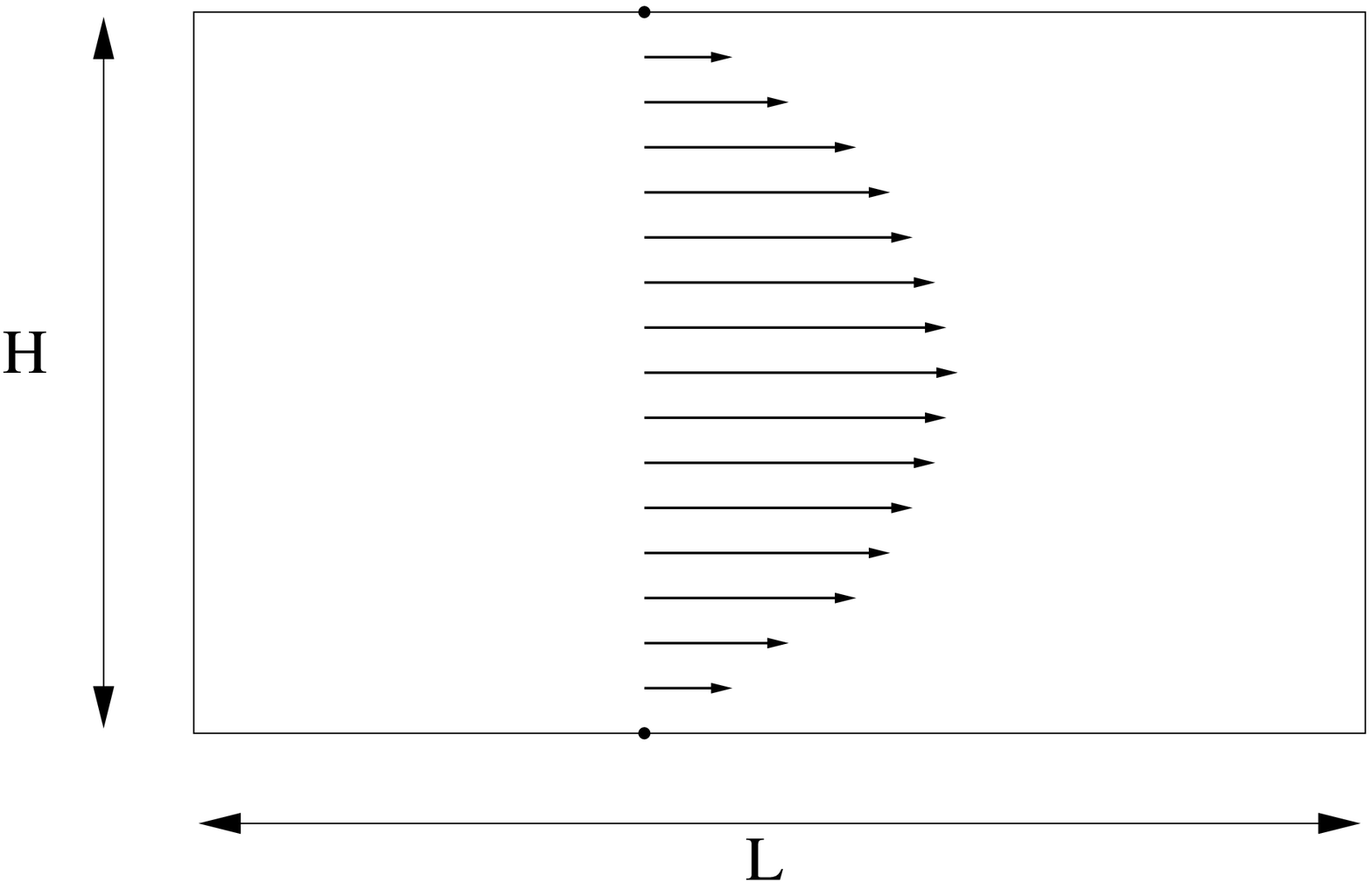}} 
\qquad    \qquad   \qquad 
   {\includegraphics [ width=.3  \textwidth  ]{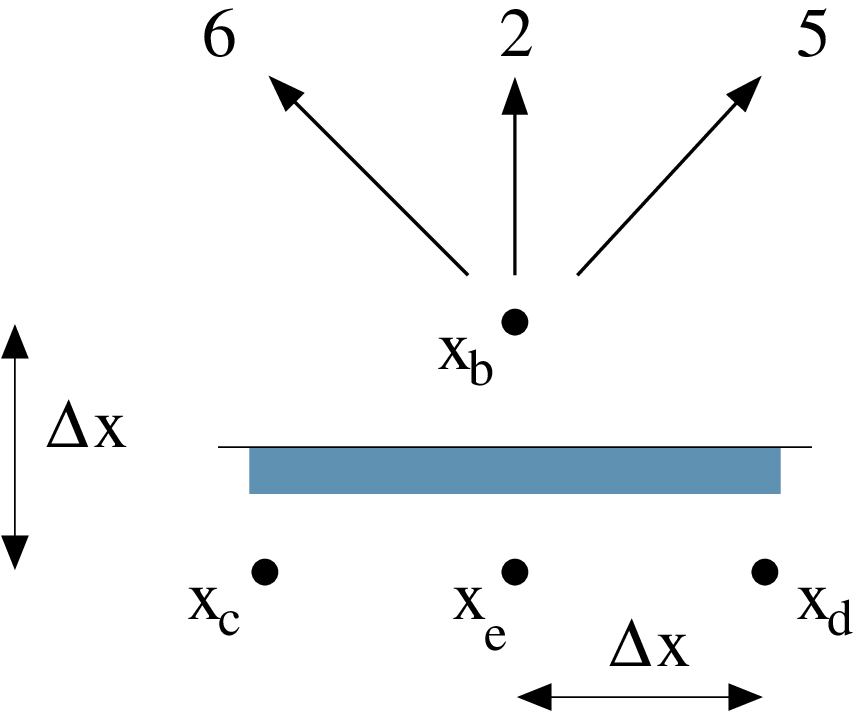}} }
\caption{Domain $\Omega \, = \, ]0,{\rm L}[  \, \times \,  ]0,H[$ (left) and notations for
  the numerical treatment of a boundary vertex $x_b$ at the bottom 
 of domain  $\Omega$ (right). }
 \label{poi}
\end{figure}
%
  
\smallskip  \noindent  $\bullet$  \quad {\bf The  Poiseuille flow} \\
We introduce a two-dimensional domain
$ \Omega=]0,{\rm L}[ \times ]0,H[$   (see Figure $\ref{poi}$).
Let ${\bf{u}}(t,x,y)\equiv (u,v)$ be the velocity of fluid and $p$
the pressure solution of the ``Poiseuille'' Stokes system:
%
%
 \begin{equation} 
   -\nu \, \Delta {\bf{u}}  + \frac{1}{\rho} \, \nabla p = 0 \,, \qquad 
  \textrm {div}  \, {\bf{u}} = 0    \qquad  {\rm in} \,\, \Omega  
 \label{st1}
\end{equation}
with the following boundary conditions:
 \begin{equation}  
p(0,y) =  - p( {\rm L} ,y)  =\delta p    \,\,\,  {\rm for} \,\,\, 0 \leq y \leq H 
\,, \quad  {\bf{u}}(x,0) =  {\bf{u}}(x,H) = 0  \,\,\,  {\rm for} \,\,\,  0 \leq x \leq 
 {\rm L}    \, . 
\label{BC1}
\end{equation} 
The solution of the above problem is classical: 
$ {\bf u}(x,\,y) \equiv  ( u(y) =  K \,y \, (H-y),  \,$  $  v =  0 ) \,, $ 
$  p(x,y)   =  2  \, \rho \, \nu \, K \, x + P_{0} \,, $ 
where 
$P_{0}$ is a  given  constant.   
We note that the problem $(\ref{st1})$ $(\ref{BC1})$ is equivalent to a
flow resulting from the action of a constant external force ${\bf{F}}$ between two
 arallel walls with  periodic boundary conditions on the inlet
and the outlet of the channel ({\it{i.e.}} in the O$x$ direction). So the
problem $(\ref{st1})$  $(\ref{BC1})$ becomes:
 \begin{equation}   
 -\nu \, \Delta {\bf{u}}  =  {\bf{F}} \,,  \qquad   
{ \bf u}  (x,0)={\bf u} (x,H)  = 0,
\label{PP3}
\end{equation}  
where ${\bf{F}}=(F_x,0)$ is the external force. The solution is given by 
$ \,  
{\bf{u}}(x,y)=  \big( u(y)=\frac{F_x}{2 \, \nu} \, y \, (H-y) \,, \,\, v =  0 \big) $.

\bigskip  \noindent $\bullet$   \quad  {\bf{A first lattice Boltzmann  scheme}} 

\noindent We use the D$2$Q$9$ lattice Boltzmann  scheme to model the Poiseuille flow
described by equation $(\ref{PP3})$. We use the equilibrium (\ref {eq-d2q9}) 
with $ \, s_{7}  = s_{8}  =  
\left(\frac{1}{2}+\frac{3\nu}{\lambda^{2}\Delta t}\right)^{-1}  \, $  
to have the exact  viscosity $\nu$  present in  equation ($\ref{PP3}$).
The implementation of the lattice Boltzmann scheme is conducted as follows
for an arbitrary mesh vertex $ \, x \, $ of the lattice.

At initial time step $t=t_{0} $ we set the vectors
$f(x,\, t_{0})=0$. 
For any given time $t$, we first determine the moments $m_k$ 
using the relation $ \, m \equiv M \, f . \,$ 
Then we
 change velocity $j_{x}$  before the collision step 
by adding a half of the external 
force $F_x$: $ \wt {\textrm \j}_{x}=j_{x}+\frac{\Delta t}{2} \, F_x $.
Thus  the macroscopic moments (density and velocity) are
evaluated.  
Then we  perform the collision step in moments space 
according to relation (\ref{moments-hors-eq}) 
and we add  half of the external force $F_x$ to the conserved velocity 
$ \wt {\textrm \j}_{x}$: 
$ j_{x} (t+\Delta t)  =   \wt {\textrm \j}_{x} + \frac{\Delta t}{2} \, F_x$.
 Using the matrix $M^{-1}$ we compute the particle distributions 
      $f_{\alpha}^{*}(x, \, t )$.
 We perform advection through a relation analogous to (\ref{dubois}) and 
 we obtain the vector  
$f_{\alpha}( x  + v_{\alpha} \, \Delta t  ,\,   t +\Delta t )$
 for $0\leq \alpha \leq 8$,
if $  x + v_{\alpha} \, \Delta t $ 
is a lattice node. 
 For a boundary node 
as $ \, x_b \,$ of Figure 3, 
we use (with the usual numbering of the degrees 
of freedom for D2Q9 scheme  \cite{bibl:La00})
the following bounce-back  boundary condition   
 \begin{equation}  \left\{ \begin{array}{rcccccl} 
  f_2 (x_b, \, t+\Delta t) &=&  f_4 (x_e, \, t+\Delta t) &=&  f_4^*  (x_b, \, t )  \\
  f_5 (x_b, \, t+\Delta t) &=&  f_7 (x_c, \, t+\Delta t) &=&  f_7^*  (x_b, \, t )  \\
  f_6 (x_b, \, t+\Delta t) &=&  f_8 (x_d, \, t+\Delta t) &=&  f_8^*  (x_b, \, t ) \,.   
\label{cl-d2q9}
 \end{array} \right.  \end{equation}
%
Periodic boundary conditions are considered in the longitudinal  direction 
for abscissae equal to $0$ and L. 
We  repeat those steps until convergence to a  steady state.

\smallskip  \noindent  {\bf Proposition 3. \quad   Order three for bounce-back }

\noindent 
For the D2Q9  lattice Boltzmann scheme   
 (\ref{dubois})  (\ref {moments-hors-eq})   (\ref {eq-d2q9} )  
  the bounce-back numerical boundary  condition  (\ref{cl-d2q9}) 
is of order $3$  at location $\Delta q = \frac{\Delta x}{2}$ 
for the Dirichlet boundary condition $ \, {\bf u} = 0 \,$ if 
and only if  $ \, \sigma_{5} \, \sigma_{8}=\frac{3}{8}. $

%

\smallskip \noindent   {\bf Proof of Proposition 3.}

\noindent    
We calculate the defects of conservation $ \theta_k $ defined  by  
(\ref{defect}) for $  k > 3 $:
 \begin{equation*}  
 \theta_{3} = \partial_t (\alpha \rho) 
\,=\, -\alpha\, {\rm{div}} j + {\rm{O}}(\Delta x^{2}) \,, \qquad  
 \theta_{4}  = \partial_t (\beta \rho) - {\rm{div}}j  
\,=\,  - (\beta+1) \, {\rm{div}}j + {\rm{O}}(\Delta x^{2})  \,, 
 \end{equation*}   
 \begin{equation*}  
\theta_{5}  =  -\frac{\partial_{t}j_{x}}{\lambda}
+\frac{\lambda(\alpha+\beta)}{3} \, \partial_{x}\rho
\,=\,   \frac{\lambda}{6} \, 
 (4+3\alpha+2\beta) \, \partial_{x}\rho+{\rm{O}}(\Delta x) \,, 
 \end{equation*}   
 \begin{equation*} 
\theta_{6} = -\frac{\partial_{t}j_{y}}{\lambda} 
+ \frac{\lambda(\alpha+\beta)}{3} \, \partial_{y}\rho  \,=\,  \frac{\lambda}{6}
\,  (4+3\alpha+2\beta) \, \partial_{y}\rho+{\rm{O}}(\Delta x) \,,    
 \end{equation*}   
 \begin{equation*} 
\theta_{7} = \frac{2}{3} \, (\partial_{x}j_{x}-\partial_{y}j_{y}) \, \qquad 
\theta_{8}=  \frac{1}{3} \, (\partial_{y}j_{x}+\partial_{x}j_{y}) \,.
 \end{equation*}   
%
Nonequilibrium moments at second order are given by the expansion (\ref{order2})   
(justified in  \cite{bibl:Du09, bibl:DLT08}). Then we have:
%
%
\begin{equation*}
m_{3}^{*} = \alpha \rho+\Delta t  \,
\, \big(\frac{1}{2}-\sigma_{3} \big) \, \big[-\alpha \, {\rm{div}}j+\Delta
			  t\frac{\lambda^{2}}{6}(\sigma_{3}\alpha(4+\alpha)
+\sigma_{5}(4+3\alpha+2\beta))\Delta \rho   \big]
+{\rm{O}}(\Delta x^{3}),
\end{equation*}
\begin{equation*}
m_{4}^{*} = \beta \rho+\Delta t \, 
(\frac{1}{2}-\sigma_{4})  \big[ -(\beta+1){\rm{div}}j +  \Delta t
\frac{\lambda^{2}}{6}(\sigma_{4}(\beta+1)(4+\alpha)+\sigma_{5}(4+3\alpha+2\beta))\Delta \rho   \big]
+{\rm{O}}(\Delta x^{3}),  
\end{equation*}
\begin{equation*}  \begin{array}{l} \displaystyle 
m_{5}^{*}  = -\frac{j_{x}}{\lambda}+\Delta t
 \big(\frac{1}{2}-\sigma_{5}\big) \, \Big[\lambda\frac{(4+3\alpha+2\beta)}{6}\partial_{x}\rho
+\Delta t \frac{\lambda}{3}  \Big(\frac{(4+3\alpha+2\beta)}{2} 
\, \sigma_5 \, \partial_x   {\rm{div}}j  
\\  [2mm]  \displaystyle \,  
+ \alpha\sigma_{3} \, \partial_{x}{\rm{div}}j+(\beta+1)  \sigma_{4} \, \partial_{x}{\rm{div}}j+2
  \sigma_{8} \partial_{x}(\partial_x  j_{x}-\partial_{y}j_{y})-\sigma_{8}\partial_{y}
  (\partial_{y}j_{x}+\partial_{x}j_{y})  \Big) \Big]+{\rm{O}}(\Delta x^{3}),
 \end{array} \end{equation*}
\begin{equation*}  \begin{array}{l} \displaystyle  
m_{6}^{*}  =  -\frac{j_{y}}{\lambda}+\Delta t 
  \big(\frac{1}{2}-\sigma_{5} \big)\, \Big[\lambda\frac{(4+3\alpha
+2\beta)}{6}\partial_{y}\rho+\Delta t \frac{\lambda}{3}  \Big(
\frac{(4+3\alpha+2\beta)}{2}  
\, \sigma_5 \, \partial_y   {\rm{div}}j  
\\  [2mm]  \displaystyle \quad 
 + \alpha \sigma_{3}\partial_{y}{\rm{div}}j+(\beta+1)   \sigma_{4}
\partial_{y}{\rm{div}}j-2  \sigma_{8} \partial_{y}(\partial_x   j_{x}
-\partial_{y}j_{y })-\sigma_{8}\partial_{x}   (\partial_{y}j_{x}+\partial_{x}j_{y}) \Big)   \Big]
+{\rm{O}}(\Delta x^{3}),
 \end{array} \end{equation*}
\begin{equation*}  \begin{array}{l} \displaystyle  
m_{7}^{*} =  \Delta t \, 
   \big(\frac{1}{2}-\sigma_{8}  \big)  \Big[\frac{2}{3}   (\partial_{x}j_{x}-\partial_{y}j_{y}) 
\\  [2mm]   \displaystyle  \qquad  \qquad   \qquad  \qquad    \qquad   
+ \, \Delta t \, \frac{\lambda^{2}}{9} \Big(\sigma_{8} (4+\alpha)  
+ \frac{ \sigma_{5}}{2} \, (4+3\alpha+2\beta)  \Big)    (\partial^{2}_{x}\rho-\partial_y^{2}\rho)
                           \Big ]+{\rm{O}}(\Delta x^{3}), 
 \end{array} \end{equation*} 
\begin{equation*}  
m_{8}^{*}  =  \Delta t \, 
   \big( \frac{1}{2}-\sigma_{8} \big)
   \Big[\frac{1}{3}(\partial_{y}j_{x}+\partial_{x}j_{y}) 
+   \Delta t \frac{\lambda^{2}}{9}
 (\sigma_{8} (4+\alpha)-\sigma_{5} (4+3\alpha+2\beta)) \partial_{xy}\rho  \Big]
+{\rm{O}}(\Delta x^{3}) \, . 
\end{equation*} 
%
We have $  j_y=0  ,\, $ $ j_{x}=j_{x0} + y \partial_{y}j_{x} + \frac{y^{2}}{2}
\partial^{2}_{y}j_{x} \, $ and $ \, \rho={\rm{constant}} $.   
We evaluate  the non conserved moments 
%
$m_{k}^{*}(x_b)$ 
and   add   $m_{1}(x_b) =j_{x}(x_b) = j_x(x_i)
+\frac{\Delta x^{2}}{3} \sigma_{8} \, \partial_{y}^{2} j_{x}(x_i) +{\rm{O}}(\Delta x^{3})  $. 
We compute moments $m_{k}^{*}(x_e)$,   $m_{k}^{*}(x_c)$ 
and  $m_{k}^{*}(x_d)$  at the ``external  nodes''  depicted in Figure~\ref{poi}. 
 Using the matrix $M^{-1}$ we evaluate  $f_k^*(x_b)$,   $f_k^*(x_c)$,   $f_k^*(x_d)$
and    $f_k^*(x_e)$. Finally we obtain  
%
\begin{equation}   
f_{5}^*(x_c) - f_{7}^*(x_b)  =  \frac{1}{6} j_{x} (x_i) + 
  \frac{ \Delta x^{2}}{144}\left(8\sigma_{5}\sigma_{8}-3\right)\frac{\partial^{2}
   j_{x}}{\partial y^{2}}(x_{i}) + {\rm{O}} (\Delta x^3) \,    
\label{magic-d2q9}
\end{equation}  
and  similar relations  for  $ \,  f_{2}^*(x_e) - f_{4}^*(x_b) $ and 
$ \, f_{6}^*(x_d) - f_{8}^*(x_b) $. 
The conclusion is clear: when the left hand side of (\ref{magic-d2q9}) is identically null
due to the boundary condition (\ref{cl-d2q9}), the momentum  $ j_{x} (x_i) $ 
on the surface located at $ \, \Delta q = \frac{1}{2} \Delta x \, $ is  null 
``up to third order accuracy'' as defined  in (\ref{def-ordre-cl})
if and only if the relation $ \, 8\sigma_{5} \, \sigma_{8} - 3 = 0 \, $ occurs.   
\hfill $\Box$   

\smallskip  \noindent $\bullet$ \quad 
We remark that if 
  we apply the body force following the algorithm 
({\it i})~$ \, m  = M \, f ,  \,\, $ 
({\it ii})~collision, $ \,\, $ 
({\it iii})~$ \, f  = M^{-1} \, m  , \,\, $  
({\it iv})~apply  the body force following the precise relations 
for  transformation of particle distribution 
$ \, f \longrightarrow \wt{f} $:  
$  \wt{f}_{1} =  f_{1}+\frac{F_x}{3\lambda} , \,   $  
$  \wt{f}_{2} =  f_{2}  , \, $  
$  \wt{f}_{3} =  f_{3}-\frac{F_x}{3\lambda} , \, $  
$  \wt{f}_{4} =  f_{4}  , \, $  
$  \wt{f}_{5} =  f_{5} + \frac{F_x}{12\lambda} , \, $  
$  \wt{f}_{6} =  f_{6} - \frac{F_x}{12\lambda} , \,$  
$  \wt{f}_{7} =  f_{7} - \frac{F_x}{12\lambda}, \,  $  
$  \wt{f}_{8} =  f_{8} + \frac{F_x}{12\lambda} \,,  $   
which are  equivalent in moments space to 
$ \,\,   \wt  {\textrm \j}_{x}   \,=\, j_{x}+F_x ,\, $  
$ \,\,    \wt{m}_{5}  \,=\, m_{5}-\frac{F_x}{\lambda} \, $   and 
$\, \wt{m}_{k}  \,=\, m_{k} \,\, $ for the other moments, 
the solid wall for the Poiseuille problem  is
 ``numerically  located''  at $\Delta  q=\frac{\Delta x}{2}$ 
up to third order accuracy  if  the  relation 
$ \, \sigma_{5} \, \sigma_{8} = \frac{3}{16} \, $ 
is satisfied between the relaxation parameters, 
as proposed by  Ginzburg and d'Humi\`eres~\cite{bibl:Gi05, bibl:GdH03, bibl:GVdH08}.

\bigskip \noindent $\bullet$  \quad  {\bf{A second  lattice Boltzmann  scheme}} \\
We can also model the Poiseuille flow described by  $(\ref{st1})$  $(\ref{BC1})$ 
 with the explicit introduction of a pressure gradient 
$ \, \delta p .$  So  the scheme (\ref{moments-hors-eq}) (\ref{eq-d2q9})
has the same steps as the preceeding scheme with 
      $ F_x \equiv 0$ and the wall boundary conditions 
are still  given by    (\ref{cl-d2q9}). 
We consider  the boundary condition for nodes 
$ \, x \equiv (k \, \Delta x , \, \ell  \,  \Delta x ) \, $ 
at the entrance ($k=1$) and at the output  ($k=N_x$) 
 as follows:
\begin{equation}  \left\{ \begin{array}{rcl}
f_1(1, \, \ell )& = & \displaystyle -f_3(0, \, \ell) + \frac{1}{18} \, (4-\alpha-2\beta) 
\, \delta \rho \,,\\
f_5(1, \, \ell )& =&  \displaystyle  -f_7(0, \, \ell-1) + \frac{1}{18} \,
(4-\alpha-2\beta)  \,  \delta \rho   \,,\\
f_8(1, \, \ell )& =&  \displaystyle-f_6(0, \, \ell +1) + \frac{1}{18} \, (4-\alpha-2\beta) 
\, \delta \rho  \,,\\ 
f_3(N_x, \, \ell )& =&  \displaystyle  -f_1(N_x+1, \, \ell)  
-   \frac{1}{18} \, (4-\alpha-2\beta)  \, \delta \rho \,, \\ 
\displaystyle  f_6(N_x, \, \ell )& =&  \displaystyle -f_8(N_x+1, 
\, \ell-1) -   \frac{1}{18} \, (4-\alpha-2\beta)  \, \delta \rho   \,,\\
f_7(N_x, \, \ell )& =&  \displaystyle-f_5(N_x+1, \, \ell +1) 
-   \frac{1}{18} \, (4-\alpha-2\beta)  \, \delta \rho  \,, 
\label{periodic-d2q9} \end{array} \right.  \end{equation}
%
with $\,\delta \rho = \, \delta p/c_s^2 \, $ the density drop
corresponding to the pressure  step  considered in (\ref{BC1}), ($\ c_s\ $ is the
speed of sound) 
and ($\alpha$,  $\beta$) parameters for equilibrium 
 introduced at the relation  (\ref{eq-d2q9}). 
Note that these expressions may be called ``anti bounce-back'' with an
imposed scalar quantity (similar to what is used when the lattice 
Boltzmann scheme  is set to simulate diffusion problems).

\smallskip  \noindent  {\bf Proposition 4. \quad   Order three for bounce-back } 

\noindent  
For the D2Q9  lattice Boltzmann scheme   
 (\ref{dubois})  (\ref {moments-hors-eq})   (\ref {eq-d2q9} )   (\ref {periodic-d2q9} ),  
  the bounce-back numerical boundary  condition  at the wall (\ref{cl-d2q9}) 
is of order $3$  at location $\Delta q = \frac{\Delta x}{2}$ 
for the Dirichlet boundary condition $ \, {\bf u} = 0 \,$ 
if and only if   $\sigma_{5} \, \sigma_{8}=-\frac{3}{8} \, \frac{\alpha+4}{\alpha+2\beta-4}.$ 
%

\smallskip \noindent   {\bf Proof of Proposition 4.}

\noindent     
In this case we perform the same proof as for proposition $3$, we take $F_x=0$
and the exact solution is given by a linear longitudinal  profile for density 
and a parabolic transverse
profile for longitudinal momentum. the algebra then follows what is presented  for
proposition~3. 
\hfill $\Box$ 
%

\smallskip \noindent $\bullet $ \quad 
 We then perform simulations of the two situations discussed above.
For this we consider a domain of size $N_x=100, N_y=21$ and analyze the flow
in the steady state. For several values of $\sigma_5$ and of $\sigma_8$, we 
determine a parabola by best fit with the velocity profile in the middle section
of the channel. We verified that the domain was long enough
in order to reduce to a negligible level the errors due to
mismatch in the end boundary conditions for links
that intersect both a solid boundary (imposed flux) 
and the input boundary (imposed pressure), that would require a
more sophisticated treatment.

\begin{figure}[!h]      
\centerline {\includegraphics [ width=.5 \textwidth, angle=-90]{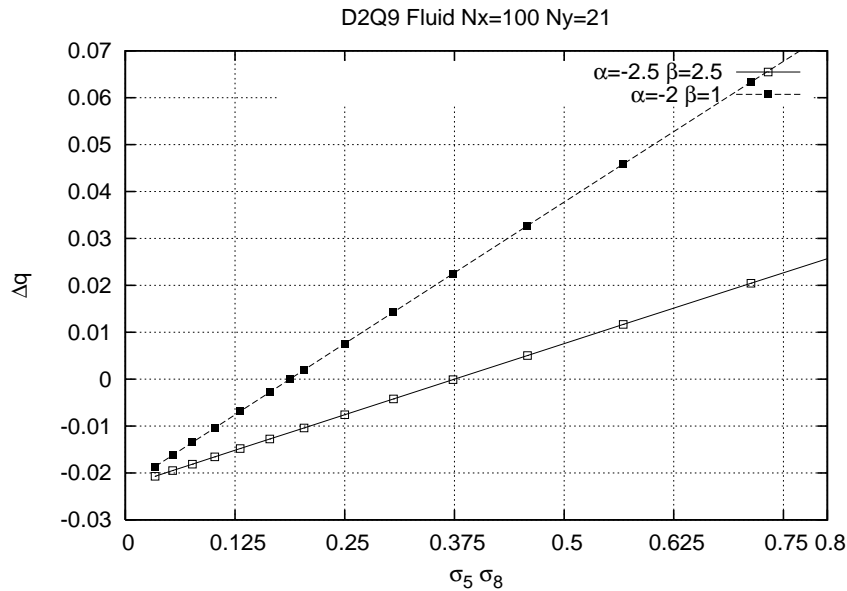}}
\caption{Product $\sigma_{5}\sigma_{8}$ {\it{versus}} solid wall  location
 $\Delta q$ with $N_{y}=21$ for the second lattice Boltzmann  boundary  scheme 
(\ref{periodic-d2q9}). 
Parameters  ${\alpha}=-2$, $\beta=1$ with 
 $ \blacksquare$ symbol,  parameters  ${\alpha}=-2.5$, $\beta=2.5$ 
 with  $ \Box$  symbol. The boundary is experimentally located at 
$\Delta q = {{\Delta x}\over{2}}$ for $\sigma_5 \, \sigma_8 = {{3}\over{16}}$ in the first
 case and $\sigma_5 \, \sigma_8 = {{3}\over{8}}$ in the second, as suggested in 
 Proposition 4.} 
\label{d2q9g}
\end{figure}

As in relation (\ref{def-ordre-cl}), we define $\Delta q$ 
as the experimental point 
where the parabola goes through zero.
The results (Figure~\ref{d2q9g})    depend only upon the product $\sigma_5 \sigma_8$ 
and are coherent with the theoretical results  established in propositions 3 and 4. For 
$\, \alpha = -2 $ and  $ \beta = 1 $,  the superconvergent accuracy is obtained 
``experimentally exactly'' at the  boundary for 
$ \,\sigma_5 \, \sigma_8 = \frac{3}{16}$.   When $\, \alpha = -\beta = -\frac{5}{2}$,
  the same observation occurs for 
 $ \,\sigma_5 \, \sigma_8 = \frac{3}{8}$.

\bigskip  \noindent {\bf \large 5 $\,\,$  Conclusion}  

\noindent 
The ``magic'' parameters introduced by Ginzburg and Adler
\cite{bibl:GA94} which allow to increase the accuracy of lattice
Boltzmann simulations in the presence of solid boundaries have been
considered for a few simple situations. We have shown that they depend
upon the choice of moments and of their equilibrium values. In addition
they depend upon the way the flow is driven. The analysis requires the
determination of the non equilibrium moments up to second order
accuracy and this explicitation is obtained 
in the framework of the Taylor expansion method.
Note that the same results could be obtained with the
Chapman-Enskog procedure. The work described here can easily be extended
to more complicated lattice Boltzmann schemes for boundaries parallel to
one of the velocities of the model. In all cases that we considered,
the results can be
expressed in terms of products of the type $\sigma_i \sigma_j$, where
$\sigma_i$ corresponds to the relevant transport coefficient
(diffusivity or shear viscosity) and $\sigma_j$ to other moments of
opposite symmetry ({\it{i.e.}} odd order moments of {$f$}, ``energy flux'' and
higher order terms of the same symmetry for models with a large enough
number of velocities),  and thus the ``magic'' conditions are the same
as those presented in the comprehensive paper of
Ginzburg, Verhaeghe and d'Humi\`eres~\cite {bibl:GVdH08}.  
They are also valid for special BGK situations that  we have in addition to the
``magic'' conditions, $\sigma_i=\sigma_j$. 
 
\bigskip  \noindent {\bf \large Acknowledgments}  

\noindent  
The referees conveyed  to the authors very interesting  remarks that have
been incorporated into  the present edition of the article.



\begin{thebibliography}{}

\end{thebibliography}


\bigskip  \noindent {\bf \large  References } 
 

\begin{thebibliography}{99} 



\bibitem{bibl:AO09} {\vskip -.3 cm}
P. Asinari, T. Ohwada.  Connection between kinetic methods for fluid-dynamic 
equations and macroscopic finite-difference schemes, 
 {\it{Computers and Mathematics with Applications}},   
 {\bf{58}}, p.~841-861,  doi:10.1016/j.camwa.2009.02.009, 2009. 
\bibitem{bibl:Ch90} {\vskip -.3 cm}
S.C. Chang. A critical analysis of the modified equation technique of Warming and
Hyett,  {\it Journal of Computational Physics},  {\bf{86}}, p.~107-126, 1990. 
\bibitem{bibl:Du07}  {\vskip -.3 cm}
F. Dubois. Une introduction au sch\'ema de
	      Boltzmann sur r\'eseau, {\textit{ESAIM: Proceedings}},
	      {\bf{18}}, p.~181--215, 2007.
\bibitem{bibl:Du08}  {\vskip -.3 cm}
F. Dubois. Equivalent partial differential equations
	of a lattice Boltzmann scheme, {\it{Computers and Mathematics
	with Applications}}, {\bf{55}}, p.~1141--1149, 2008.  
\bibitem{bibl:Du09}  {\vskip -.3 cm}
F. Dubois.  Third order equivalent equation  of  lattice 
Boltzmann scheme,    {\it Discrete and Continuous Dynamical Systems-Series A},  
  {\bf 23},  p.~221-248, 2009. 
\bibitem{bibl:DLT08}  {\vskip -.3 cm}
F. Dubois, P. Lallemand, M. M. Tekitek. 
	       Using the Lattice Boltzmann Scheme for Anisotropic
	       Diffusion Problems,   {\textit{Finite Volumes for complex
	Applications V}, R. Eymard, J.M. H\'erard (Eds)},  p.~351--358, Wiley,   2008.

\bibitem{bibl:FHH87}  {\vskip -.3 cm}
U. Frisch, D. d'Humi{\`e}res, B. Hasslacher,
	P. Lallemand, Y. Pomeau, J.-P. Rivet.  {Lattice gas hydrodynamics
	in two and three dimensions}, {\it{Complex Systems}}, {\bf{1}},
	p.~{649--707},  1987.
\bibitem{bibl:Gi05}  {\vskip -.3 cm}
I. Ginzburg. Generic boundary conditions
	      for lattice Boltzmann models and their application to
	      advection and anisotropic dispersion equations,
	      {\textit{Advances in Water Resources}}, {\bf{28}},
	      p.~1196-1216, 2005.
\bibitem{bibl:GA94}  {\vskip -.3 cm}
I. Ginzbourg, P.M. Adler. Boundary flow condition
	analysis for three-dimensional lattice Boltzmann model,
	{\it{J. Phys. II France}}, {\bf{4}}, p.~191--214, 1994.
\bibitem{bibl:GD95}  {\vskip -.3 cm}
I. Ginzburg, D. d'Humi\`eres. Second order boundary method for Lattice 
Boltzmann model, 
{\it Journal of  Statistical  Physics},   {\bf 84}, p.~927-971, 1995. 
\bibitem{bibl:GdH03}  {\vskip -.3 cm}
I. Ginzburg, D. d'Humi\`eres.  Multireflection
	boundary conditions for lattice Boltzmann models,
	{\it{Phys. Rev. E}}, {\bf{68}}, p.~66614--66644, 2003.
\bibitem{bibl:GVdH08}  {\vskip -.3 cm}
I. Ginzburg, F. Verhaeghe and D. d'Humi\`eres. Two-relaxation-time lattice
	Boltzmann scheme: About parametrization, velocity, pressure and
	mixed boundary conditions, {\it{Communications in Computational
	Physics}}, {\bf{3}}, p.~519-581, 2008. 
\bibitem{GSS86}  {\vskip -.3 cm}
D. Griffiths, J. Sanz-Serna. On the scope of the method of modified equations, 
{\it SIAM Journal on  Scientific and Statistical  Computing}, {\bf 7}, p.~994-1008, 1986.
\bibitem{bibl:HJ89} {\vskip -.3 cm}
F. Higuera, J. Jim\'enez. Boltzmann approach to lattice gas simulations, 
 {\it Europhysics Letters}, {\bf 9},  
  p.~663-668, 1989.  
\bibitem{bibl:HSB89} {\vskip -.3 cm}
F. Higuera, S. Succi and R. Benzi. Lattice gas dynamics with
  enhanced collisions,  
{\it Europhysics Letters},  {\bf 9},  n$^{\rm o}$~4,  p.~345-349, 1989.
\bibitem{bibl:dh92}  {\vskip -.3 cm}
D. d'Humi\`eres. Generalized Lattice-Boltzmann Equations, 
in {\it Rarefied Gas Dynamics: Theory
and Simulations}, {\it AIAA Progress in Astronautics and Astronautics},
  {\bf 159},  p.~450-458, 1992. 
\bibitem{bibl:JKL05}    {\vskip -.3 cm}
 M. Junk, A. Klar, and L.-S. Luo. Asymptotic analysis of the lattice 
Boltzmann equation, {\it 
Journal of Computational Physics},  {\bf 210},  p.~676-704, 2005. 
\bibitem{JY09a}    {\vskip -.3 cm}
 M. Junk, Z. Yang. Convergence of Lattice Boltzmann Methods 
for Navier-Stokes Flows in Periodic and Bounded Domains, 
{\it Numerische Mathematik},  {\bf 112},  p.~65-87, 2009. 
 \bibitem{jy09b} {\vskip -.3 cm}
 M. Junk, W.-A. Yong.  Weighted L$^2$ Stability of the Lattice Boltzmann 
Method, {\it SIAM Journal on Numerical Analysis},  {\bf 47}, p.~1651-1665, 2009.
\bibitem{bibl:KGSB98}   {\vskip -.3 cm}
I.V. Karlin, A.N. Gorban,  S. Succi and V. Boffi. 
Maximum Entropy Principle for Lattice Kinetic Equations, 
 {\it Physical  Review Letters}, {\bf 81},    p.~6-9, 1998. 
\bibitem{bibl:La00}  {\vskip -.3 cm}
P. Lallemand, L. Luo. Theory of the lattice
 Boltzmann method: Dispersion, dissipation, isotropy, Galilean
 invariance, and stability, {\em{Physical Review E}}, {\bf{61}},
	p.~6546--6562, 2000. 
\bibitem{bibl:LP74}  {\vskip -.3 cm}
A. Lerat, R. Peyret. Noncentered Schemes and Shock Propagation
  Problems,  {\it  Computers and Fluids},  {\bf 2}, p.~35-52, 1974. 
\bibitem{bibl:qian} {\vskip -.3 cm}
Y.H. Qian, D. d'Humi\`eres, P. Lallemand.
Lattice BGK models for Navier-Stokes equation, 
{\em{Europhys. Lett.}}, {\bf{17}}, p.~479--484, 1992.
\bibitem{bibl:VR99} {\vskip -.3 cm}
F.R. Villatoro, J.I. Ramos. 
On the method of modified equations. V: 
Asymptotic analysis of and direct-correction and asymptotic 
successive-correction techniques for the implicit midpoint method, 
{\it Applied Mathematics and Computation},  {\bf{103}}, p.~241-285, 1999.  
\bibitem{bibl:WH74}   {\vskip -.3 cm}
R.F. Warming,  B.J. Hyett. The modified equation approach 
to the stability and accuracy analysis of finite difference methods, 
{\it Journal of Computational Physics},  {\bf 14}, p.~159-179, 1975. 
 





\end{thebibliography}
\end{document}